\documentclass[11pt,reqno]{amsart}
\usepackage{etex}
\usepackage{dsfont,amsmath,amsfonts,amscd,amssymb,graphicx,mathrsfs,eufrak,upgreek}
\usepackage[dvips,all,arc,curve,color,frame]{xy}
\usepackage[usenames]{color}
\usepackage{setspace}
\usepackage{array,multirow,booktabs,longtable}

\usepackage[final]{optional}

\usepackage{versions}
\includeversion{arXiv}

\newcommand{\tocspace}{0.1ex}
\let\oldtocsection=\tocsection
\let\oldtocsubsection=\tocsubsection
\let\oldtocsubsubsection=\tocsubsubsection
\renewcommand{\tocsection}[3]{\hspace{0em}\oldtocsection{#1}{#2}{#3}\vspace{\tocspace}}
\renewcommand{\tocsubsection}[3]{ \hspace{1em} \oldtocsubsection{#1}{\small{#2}}{\small{#3}}\vspace{\tocspace} }
\renewcommand{\tocsubsubsection}[3]{\hspace{2em}\oldtocsubsubsection{#1}{\small{#2}}{\small{#3}}}

\newcommand{\marginparstretch}{0.6}
\let\oldmarginpar\marginpar
\renewcommand\marginpar[1]{\-\oldmarginpar[\framebox{\setstretch{\marginparstretch}\begin{minipage}{\marginparwidth}{\raggedleft\tiny #1}\end{minipage}}]{\framebox{\setstretch{\marginparstretch}\begin{minipage}{\marginparwidth}{\raggedright\tiny #1}\end{minipage}}}}

\usepackage[colorlinks]{hyperref}
\usepackage{tikz,tikz-cd,mathrsfs}
\usetikzlibrary{arrows,decorations.pathmorphing,decorations.pathreplacing,positioning,shapes.geometric,shapes.misc,decorations.markings,decorations.fractals,calc,patterns}

\tikzset{
        Point/.style={circle,draw=black,circle,fill=black,inner sep=0pt, minimum size=2pt},
        DynkinBlack/.style={circle,draw=black,circle,fill=black,inner sep=0pt, minimum size=4pt},
         DynkinWhite/.style={circle,draw=black,circle,fill=white,inner sep=0pt, minimum size=4pt},
        cvertex/.style={circle,draw=black,inner sep=1pt,outer sep=3pt},
        vertex/.style={circle,fill=black,inner sep=1pt,outer sep=3pt},
        star/.style={circle,fill=yellow,inner sep=0.75pt,outer sep=0.75pt},
        tvertex/.style={inner sep=1pt,font=\scriptsize},
        gap/.style={inner sep=0.5pt,fill=white}}

\tikzstyle{mybox} = [draw=black, fill=blue!10, very thick,
    rectangle, rounded corners, inner sep=10pt, inner ysep=20pt]
\tikzstyle{boxtitle} =[fill=blue!50, text=white,rectangle,rounded corners]

\newcommand{\arrow}[2][20]
 {
  \hspace{-5pt}
  \begin{tikzpicture}
   \node (A) at (0,0) {};
   \node (B) at (#1pt,0) {};
   \draw [#2] (A) -- (B);
  \end{tikzpicture}
  \hspace{-5pt}
 }

\newcommand{\birational}[1][20]{\arrow[#1]{->,dashed}}

\newtheorem{theorem}{Theorem}[section]
\newtheorem{keytheorem}{Theorem}[section]

\newtheorem{proposition}[theorem]{Proposition}

\newtheorem{definition}[theorem]{Definition}

\theoremstyle{definition} 

\newtheorem{example}[theorem]{Example}

\newtheorem{remark}[theorem]{Remark}

\theoremstyle{remark}
\newtheorem*{acknowledgements}{Acknowledgements}
\newtheorem*{outline}{Outline}

\numberwithin{equation}{section}

\newcounter{enumeratenoindentcounter}

\renewcommand{\t}[1]{\textnormal{#1}}

\def\coh{\mathop{\rm coh}\nolimits}

\def\Hom{\mathop{\rm Hom}\nolimits}

\def\End{\mathop{\rm End}\nolimits}
\def\Ext{\mathop{\rm Ext}\nolimits}

\def\Spec{\mathop{\rm Spec}\nolimits}

\def\D{\mathop{\rm{D}^{}}\nolimits}

\def\Id{\mathop{\rm{Id}}\nolimits}

\newcommand\con{\mathrm{con}}

\newcommand{\CA}{\mathrm{A}_\con}

\newcommand{\CAab}{\mathrm{A}^{\ab}_\con}
\newcommand{\AB}{\mathrm{A}}

\def\ab{\mathop{\rm ab}\nolimits}

\newcommand\art{\mathsf{Art}}

\newcommand\Sets{\mathsf{Sets}}
\newcommand\cDef{c\mathcal{D}ef}
\newcommand\Def{\mathcal{D}ef}

\newcommand{\cE}{\mathcal{E}}
\newcommand{\cF}{\mathcal{F}}

\newcommand{\cO}{\mathcal{O}}

\newcommand\surface{S}
\newcommand\surfaceRes{\widehat{S}}
\newcommand\simRes{\widehat{Y}}
\newcommand\DynkinDiag{\Upgamma}
\newcommand\chosenVertices{\DynkinDiag'}
\newcommand\TestAlg{B}
\newcommand\WeylSubgroup{W'}
\newcommand\threefold{X}
\newcommand\threefoldCon{X_\con}
\newcommand\contraction{f}
\newcommand\threefoldDef{\tilde{X}}
\newcommand\threefoldDefCon{\tilde{X}_\con}
\newcommand\hypersurface{H}
\newcommand\hypersurfaceEq{\upphi}
\newcommand\hypersurfacePRes{X_H}
\newcommand\hypersurfaceRes{\widehat{H}}
\newcommand\hypersurfaceResMap{\contraction} 
\newcommand\discriminantLocus{L}
\newcommand\pathParameter{w}
\newcommand\neighbourhood{\Updelta}
\newcommand\GV{Gopakumar--Vafa }

\newlength\tempWidth

\newcommand\reversecolon{\,:\!} 

\newcommand{\defColor}{gray!50}
\newcommand{\defBGColor}{gray!10}
\newcommand{\defGridColor}{gray}
  \def\size{5}
  \def\bulge{1/3}
\def\gridShift{-4}
\def\bend{0.1}
\def\len{0.8}
\def\dotsize{0.1}

\newcommand{\defCurveCoords}[3]{
  (#1+\bend,-#3,#2) .. controls (#1-\bend,-#3/3,#2) and (#1-\bend, #3/3,#2) .. (#1+\bend,#3,#2)
}

\newcommand{\defDrawDot}[2]{
    \draw[fill] (#1+\dotsize,\gridShift,#2+\dotsize) -- (#1+\dotsize,\gridShift,#2-\dotsize) -- (#1-\dotsize,\gridShift,#2-\dotsize) -- (#1-\dotsize,\gridShift,#2+\dotsize) -- cycle
    }

\newcommand{\defCurve}[3]{
  \draw[color=#3] \defCurveCoords{#1}{#2}{\len}
}

\newcommand{\defThreeFold}[3]{
  \draw[color=#3] (#1+0.5+\size,#2-0.3) .. controls (#1+0.5+2*\size*\bulge,#2-0.3+\size/2) and (#1+0.5-\size*\bulge,#2-0.3+\size/2) .. (#1+0.5-\size,#2-0.3) .. controls (#1+0.5-5*\size*\bulge,#2-0.3-\size/2) and (#1+0.5+4*\size*\bulge,#2-0.3-\size/2) .. (#1+0.5+\size,#2-0.3)
}

\newcommand{\defGridLines}[4]{
  \def\x{#1} \def\y{#2}
  \def\numx{#3} \def\numy{#4}
  \foreach \i in {0,...,\numy} {
    \draw[color=\defGridColor] (\x,\gridShift,\y+\i) -- +(\numx,0,0);
  };
  \foreach \i in {0,...,\numx} {
    \draw[color=\defGridColor] (\x+\i,\gridShift,\y) -- +(0,0,\numy);
  };
  \draw[dotted] (\x,0.35*\gridShift,\y) -- (\x,0.9*\gridShift,\y)
}

\newcommand{\defUniversal}[4]{
  \def\x{#1} \def\y{#2}
  \def\numx{#3} \def\numy{#4}
  \def\zero{0}
  \draw[color=\defBGColor,fill=\defBGColor] \defCurveCoords{\x}{\y}{\len} -- \defCurveCoords{\x}{\y+\numx}{-\len} -- cycle;
  \draw[color=\defBGColor,fill=\defBGColor] \defCurveCoords{\x}{\y}{\len} -- \defCurveCoords{\x+\numy}{\y}{-\len} -- cycle;
  \draw[color=\defBGColor,fill=\defBGColor] \defCurveCoords{\x+\numy}{\y}{\len} -- \defCurveCoords{\x+\numy}{\y+\numx}{-\len} -- cycle;
  \draw[color=\defBGColor,fill=\defBGColor] \defCurveCoords{\x}{\y+\numx}{\len} -- \defCurveCoords{\x+\numy}{\y+\numx}{-\len} -- cycle;
  \defUniversalCurves{#1}{#2}{#3}{#4}
  \foreach \i in {1,...,\numy} {
    \foreach \j in {1,...,\numx} {
      \defDrawDot{\x+\i}{\y+\j};
    };
  }
}

\newcommand{\defUniversalComm}[4]{
  \def\x{#1} \def\y{#2}
  \def\numx{#3} \def\numy{#4}
  \def\zero{0}
  \draw[color=\defBGColor,fill=\defBGColor] \defCurveCoords{\x}{\y}{\len} -- \defCurveCoords{\x}{\y+\numx}{-\len} -- cycle;
  \draw[color=\defBGColor,fill=\defBGColor] \defCurveCoords{\x}{\y}{\len} -- \defCurveCoords{\x+\numy}{\y}{-\len} -- cycle;
  \defUniversalCurves{#1}{#2}{#3}{#4}
}

\newcommand{\defUniversalCurves}[4]{
  \def\x{#1} \def\y{#2}
  \def\numx{#3} \def\numy{#4}
  \def\zero{0}
  \foreach \i in {0,...,#4} {
    \ifx\i\zero \defCurve{#1+\i}{#2}{black} \else \defCurve{#1+\i}{#2}{\defColor} \fi;
    \defDrawDot{#1+\i}{#2};
  };
  
  \ifx\numx\zero {} \else {
  \foreach \j in {1,...,#3} {
    \defCurve{#1}{#2+\j}{\defColor};
    \defDrawDot{#1}{#2+\j};
    }
    } \fi
}

\begin{document}

\title[Contractions of $3$-folds: deformations and invariants]{\textsc{Contractions of $3$-folds:\\ deformations and invariants}}

\author{Will Donovan}
\address{Kavli Institute for the Physics and Mathematics of the Universe (WPI), Todai Institutes for Advanced Study, The University of Tokyo, 5-1-5 Kashiwanoha, Kashiwa, Chiba, 277-8583, Japan}
\email{will.donovan@ipmu.jp}
\begin{abstract}This note discusses recent new approaches to studying flopping curves on $3$-folds. In a joint paper \cite{DW1}, the author and Michael Wemyss introduced a $3$-fold invariant, the contraction algebra, which may be associated to such curves. It characterises their geometric and homological properties in a unified manner, using the theory of noncommutative deformations. Toda has now clarified the enumerative significance of the contraction algebra for flopping curves, calculating its dimension in terms of \GV invariants~\cite{TodaGV}. Before reviewing these results, and others, I give a brief introduction to the rich geometry of flopping curves on $3$-folds, starting from the resolutions of Kleinian surface singularities.\begin{arXiv} This is based on a talk given at VBAC~2014 in Berlin.\end{arXiv}
\end{abstract}

\thanks{The author is supported by World Premier International Research Center Initiative (WPI Initiative), MEXT, Japan, and thanks the organisers of VBAC 2014 in Berlin for their kind hospitality, and for a stimulating conference.}
\maketitle
\parindent 20pt
\parskip 0pt


\section{Introduction}
\label{section intro}

In joint work \cite{DW1}, the author and Michael~Wemyss associated a certain algebra to a contraction of a smooth rational curve in a complex $3$-fold, and began an investigation of the importance of this \emph{contraction algebra} in studying the geometry of the $3$-fold. In particular, the algebra was shown to represent a certain noncommutative deformation functor associated to the curve. Further recent work extends this to general $3$-fold flopping contractions~\cite{DW3}, and to certain higher-dimensional contractions~\cite{DW2}, but for simplicity I focus in this note on flopping contractions of $3$-folds with exceptional locus a smooth rational curve.

Given a flopping contraction, the $3$-fold $X$ which is contracted  has a birational partner $X'$ called a flop. For smooth $X$, work of Bridgeland then yields an associated canonical symmetry of the derived category of $X$, known as the flop--flop functor. We use the contraction algebra to give a description of this symmetry, intrinsic to the geometry of $X$.

This same contraction algebra generalises and unifies known invariants of flopping curves, including the width invariant of Reid \cite{Pagoda} and the Dynkin type \cite{KM,Kawa}. Recent work of Toda shows that the dimension of the contraction algebra is determined by the genus zero \GV invariants associated to the contracted curve \cite{TodaGV}. In the other direction, the contraction algebra determines these \GV invariants for curves of Dynkin types~$A$ and~$D$.

For flopping curves with Dynkin type~$A$, the contraction algebra turns out to be of the simple form $\mathbb{C}[\upepsilon]/\upepsilon^d$, whereas for Dynkin types~$D$ and $E$, it is always noncommutative. This note gives a compact introduction to the geometry of flopping curves, along with simple examples of contraction algebras.

I now give preliminary statements, deferring details and definitions until Section~\ref{section threefolds}. Consider a contraction $X \to X_{\con}$ of $3$-folds which maps a curve $C\cong \mathbb{P}^1$ to a point $p$, and is an isomorphism elsewhere. We give an explicit construction of a $\mathbb{C}$-algebra~$\CA$ (Definition~\ref{Definition contraction algebra}), which we refer to as the contraction algebra, and prove the following theorem.

\begin{keytheorem}[\ref{theorem representability}]\cite{DW1} The $\mathbb{C}$-algebra $\CA$ represents a functor of noncommutative deformations associated to the curve $C$.
\end{keytheorem}

For flopping contractions, under the assumption that $X$ is projective with at worst Gorenstein terminal singularities, the following is a result of Toda.

\begin{keytheorem}[\ref{theorem.toda}] \label{theorem Toda intro} \cite{TodaGV} For a flopping curve $C$, the dimension $\dim(\CA)$ may be calculated from genus zero \GV invariants $n_j$, and the length invariant $l$, by the formula
\[
\dim(\CA) = \sum_{j=1}^l j^2 n_j.
\]
Furthermore, the dimension of the abelianisation $\dim(\CAab)$ is $n_1$.
\end{keytheorem}

Theorem~\ref{theorem Toda intro} gives a new enumerative criterion for the noncommutativity of the algebra~$\CA$, namely that it is noncommutative precisely when one of the higher \GV invariants $n_j$ for $j>1$ is non-zero.

Toda's proof of Theorem~\ref{theorem Toda intro} makes use of the following result, which gives a description of the flop--flop functor. Stating and proving this was a key motivation for the definition of $\CA$.

\begin{keytheorem}[\ref{theorem.twist}]\cite{DW1} Assuming $X$ is projective with at worst Gorenstein terminal singularities, the Bridgeland--Chen flop--flop functor \[F' \circ F \colon \D(X) \overset{\sim}{\longrightarrow} \D(X') \overset{\sim}{\longrightarrow} \D(X)\] is inverse to a symmetry $T_C$ of the derived category \,$\D(X)$ constructed using the algebra $\CA$.
The action of $T_C$ on $E \in \D(X)$ may be calculated by the following formula, where $\cE \in \coh(\cO_X \otimes_{\mathbb{C}} \CA)$ denotes the universal family for the noncommutative deformations of $\cO_C(-1)$.
\[T_C(E) \cong \operatorname{Cone}\left( \mathbf{R}\!\Hom_X(\cE,E)\otimes_{\CA}^{\bf L} E \overset{\operatorname{ev}}{\longrightarrow} E \right)\]
\end{keytheorem}

\begin{outline} In \S\ref{section surfaces}, I give some orientation on surface singularities and their resolutions and deformations, which have an elegant representation-theoretic description. I~explain in \S\ref{section threefolds} how this may be used to understand contractions of $3$-folds, by viewing them as families of contractions of surfaces. In \S\ref{section contraction algebra}, two approaches to the contraction algebra $\CA$ are given, via a tilting algebra associated to the contraction, and via noncommutative deformations. In \S\ref{section numerical invariants}, I discuss enumerative invariants associated to $\CA$. Finally in \S\ref{section derived category}, the contraction algebra is related to the flop--flop functor, and I indicate how Toda uses this relation to prove Theorem~\ref{theorem Toda intro} above. No~knowledge of the derived category is assumed until this final section.\end{outline}

\begin{acknowledgements}My first thanks go to Michael Wemyss for the great experience of working together on our paper \cite{DW1}, and its sequels~\cite{DW2,DW3}. Thanks also go to Yukinobu Toda for his interest in this work, and for sharing his insights; to Agnieszka Bodzenta, Alexei Bondal, Mikhail Kapranov, and Edward Witten for helpful conversations on the material in this note; to Joseph Karmazyn for his comments on the document; and to the anonymous referee for their valuable suggestions.

I remain grateful for the support of EPSRC grant~EP/G007632/1, and of Iain~Gordon, as well as the hospitality of the Erwin Schr\"odinger Institute, Vienna, and the Hausdorff Institute, Bonn, during the course of my work on \cite{DW1}, and for the support of the Korean Institute of Advanced Study, Seoul, and Miles Reid, while beginning work on this document.\end{acknowledgements}

\section{Resolutions of surfaces}
\label{section surfaces}

\subsection{Du Val resolutions}\label{du val} Take a surface $\surface$ with a Kleinian singularity.
Such singularities are given complete locally by quotients $\mathbb{A}^2 / G$, with $G$ a finite subgroup of $\operatorname{SL}(2,\mathbb{C})$.
These subgroups are classified by simply-laced Dynkin diagrams~$\DynkinDiag$ \cite{McKay, ReidMcKay}. The types~$A$~and~$D$ correspond to cyclic and binary dihedral groups respectively, and the types~$E_6$, $E_7$, and~$E_8$ correspond to binary Platonic groups. There exists a resolution of singularities $\surfaceRes \to \surface$ obtained by replacing the singular point with a tree of projective lines, whose intersections are determined by the diagram~$\DynkinDiag$~\cite{DuVal}. This is shown in Figure~\ref{figure minimal resolutions surfaces} for the $A_2$ and $D_4$ cases. Here the relevant groups are cyclic of order~3 and binary dihedral of order~8 respectively, with the latter being isomorphic to the quaternion group.

\newcommand{\defCurveOverlap}{0.1}

\newcommand{\defTypeADynkin}[1]{
  \node (0) at (0,0) [DynkinWhite] {};
}

\newcommand{\defTypeATwoDynkin}[1]{
  \node (0) at (0,0) [DynkinWhite] {};
  \node (1) at (1,0) [DynkinWhite] {};

  \draw [-] (0) -- (1);
}

\newcommand{\defTypeDDynkin}[1]{
  \node (-1) at (-1,0) [DynkinWhite] {};
  \node (0) at (0,0) [DynkinWhite] {};
  \node (1) at (1,0) [DynkinWhite] {};
  \node (top) at (0,1) [DynkinWhite] {};

  \draw [-] (-1) -- (0);
  \draw [-] (0) -- (1);
  \draw [-] (0) -- (top);
}

\newcommand{\defTypeDDynkinSingle}[1]{
  \node (-1) at (-1,0) [DynkinBlack] {};
  \node (0) at (0,0) [DynkinWhite] {};
  \node (1) at (1,0) [DynkinBlack] {};
  \node (top) at (0,1) [DynkinBlack] {};

  \draw [-] (-1) -- (0);
  \draw [-] (0) -- (1);
  \draw [-] (0) -- (top);
}

\newcommand{\defTypeACurves}[2]{
  \foreach \i in {1,...,#1} {
    \draw (\i-\defCurveOverlap+3/2-#1/2+#2/2,0,0) to [bend left=30] (\i+1+\defCurveOverlap+3/2-#1/2+#2/2,0,0);
  };
}

\newcommand{\defTypeDCurves}[1]{
  \defTypeACurves{#1}{0}
  \draw (2.3+3/2-#1/2,0,0+\defCurveOverlap) to [bend left=30] (2.3+3/2-#1/2,0,-1.5-\defCurveOverlap);
}

\newcommand{\defSurface}[2]{
  \draw [#2] (#1-0.5,0,0) ellipse (#1/2+1 and 1);
}

\newcommand{\defThreefold}[1]{
  \draw (#1-0.5,0,0) ellipse (#1/2+1 and 2);
}

\newcommand{\defTypeASurfaceShifted}[2]{
  \defSurface{3}{}
  \defTypeACurves{#1}{#2}
}

\newcommand{\defTypeASurface}[2]{
  \defSurface{3}{#2}
  \defTypeACurves{#1}{0}
}

\newcommand{\defTypeDSurface}[1]{
  \defSurface{3}{}
  \defTypeDCurves{#1}
}

\newcommand{\defSingSurface}[2]{
  \defSurface{3}{#2}
  \node (sing) at (3-0.5,0,0) [Point] {};
}

\newcommand{\defSurfaceScale}{0.6}
\newcommand{\defDynkinScale}{0.4}
\newcommand{\defDynkinSep}{1.7}
\begin{figure}[ht]
\begin{center}
\begin{tikzpicture}[scale=0.7] 
\node (A) at (0,3.5) {\begin{tikzpicture}[scale=\defSurfaceScale]\defTypeASurface{2}{}\end{tikzpicture}};
\node (Acon) at (8,3.5) {\begin{tikzpicture}[scale=\defSurfaceScale]
\defSingSurface{3}{}\end{tikzpicture}};
\node (ADynkin) at (0,3.5+\defDynkinSep-0.25) {\begin{tikzpicture}[scale=\defDynkinScale]
\defTypeATwoDynkin{1}\end{tikzpicture}};

\node (D) at (0,0) {\begin{tikzpicture}[scale=\defSurfaceScale]\defTypeDSurface{3}\end{tikzpicture}};
\node (Dcon) at (8,0) {\begin{tikzpicture}[scale=\defSurfaceScale]\defSingSurface{3}{}\end{tikzpicture}};
\node (DDynkin) at (0,0+\defDynkinSep) {\begin{tikzpicture}[scale=\defDynkinScale]
\defTypeDDynkin{3}\end{tikzpicture}};

\draw[->] (A) -- (Acon);
\draw[->] (D) -- (Dcon);
\node[left of=A] (labelA) at (A.west) {$A_2$};
\node[left of=D] (labelD) at (D.west) {$D_4$};
\end{tikzpicture}
\end{center}
\caption{Resolutions of Kleinian surface singularities}
\label{figure minimal resolutions surfaces}
\end{figure}
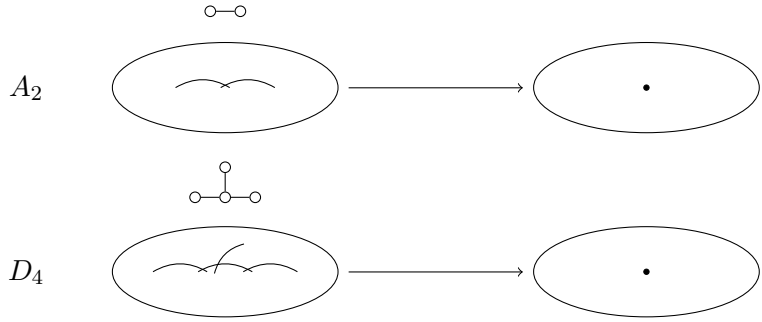

\subsection{Versal deformations} Our goal is to understand $3$-fold geometries via $1$-parameter deformations of these resolutions. Such deformations are controlled by versal deformations, which have a representation-theoretic description as follows. Associated to the diagram $\DynkinDiag$ is a complex semisimple Lie algebra, and a choice of Cartan subalgebra, denoted by $U$, carries an action of a Weyl group $W$. The singular surface $S$ has a versal deformation $Y$ parametrised by the space of $W$-orbits in $U$. After base change from this parameter space $U/W$ to $U$ itself, the fibres of the family $Y$ have a simultaneous resolution \cite{Slo} with central fibre recovering the resolution~$\surfaceRes \to \surface$. Denoting this by $\simRes$, there is a Cartesian square as follows, and the resolution $\surfaceRes \to \surface$ is recovered by taking the fibre of the top line over $0 \in U/W$.
\[
\begin{tikzcd}
\simRes \ar{r} \ar{d}  & Y \ar{d} \\
U \ar{r}  & U/W
\end{tikzcd}
\]

\subsection{Case $A_1$}\label{section A1 case} I explain this case in detail as it immediately gives the simplest $3$-fold geometry which will be of interest to us. We may take $\surface = \mathbb{A}^2 / G,$ where $G = \mathbb{Z}/2\mathbb{Z},$ acting by $\pm \Id$. Taking suitable quadratic invariants of this $G$-action, namely $x=u^2-v^2$, $y=2uv$, and $z=u^2+v^2$, the singularity $S$ may be written as a hypersurface \begin{equation*} S= \left\{ x^2 + y^2 = z^2 \right\} \subset \mathbb{C}^3,\end{equation*} and has a versal deformation \begin{equation*} Y= \left\{ x^2 + y^2 = z^2 - q\right\}  \subset \mathbb{C}^4,\end{equation*} with parameter $q$. The total space of this deformation $Y$ is smooth. 

Now in this case we have $U = \mathbb{C}$ with Weyl group $W = S_2$ acting by $\pm \Id$. Writing $q=t^2$ for a new parameter $\pm t \in \mathbb{C} / S_2$, we may view $Y$ as a family over~$U/W$. Base changing to the parameter $t\in \mathbb{C}$, a new deformation \begin{equation}\label{equation conifold} Y' = \left\{  x^2 + y^2 = z^2 - t^2\right\}\subset \mathbb{C}^4\end{equation} is obtained, with a singularity at $0$. This has a simultaneous resolution $\simRes$ which can be constructed by a blowup in one of the ideals $(x+iy,z\pm t)$, and is sketched in Figure~\ref{figure sim res A1}. This $\simRes$ is our first example of a smooth $3$-fold with a contraction. It contains a single complete curve $C$ in the fibre over $t=0$, with normal bundle $\cO_C(-1) \oplus \cO_C(-1)$. This is an example of the simplest sort of $3$-fold contractible curve, and is referred to as a $(-1,-1)$-curve.

\newcommand{\LabelDrop}{2}
\newcommand{\LabelSep}{3}
\begin{figure}[ht]
\begin{center}
\begin{tikzpicture}[scale=0.7]
\node (A) at (0,-4) {$\begin{array}{c}
\begin{tikzpicture}[scale=\defSurfaceScale]\defSurface{3}{}\end{tikzpicture} \\
\begin{tikzpicture}[scale=\defSurfaceScale]\defTypeASurface{1}{}\end{tikzpicture} \\
\begin{tikzpicture}[scale=\defSurfaceScale]\defSurface{3}{}\end{tikzpicture}
\end{array}$};
\node (Acon) at (8,-4) {$\begin{array}{c}
\begin{tikzpicture}[scale=\defSurfaceScale]\defSurface{3}{}\end{tikzpicture} \\
\begin{tikzpicture}[scale=\defSurfaceScale]\defSingSurface{3}{}\end{tikzpicture} \\
\begin{tikzpicture}[scale=\defSurfaceScale]\defSurface{3}{}\end{tikzpicture}
\end{array}$};

\draw[->] (A) -- (Acon);
\node (labelA) at (0-\LabelSep,-\LabelDrop) {$\simRes$}; 
\node (labelAcon) at (8+\LabelSep,-\LabelDrop) {$Y'$}; 
\end{tikzpicture}
\end{center}
\caption{Simultaneous resolution $\simRes$ of deformation of $A_1$ singularity}
\label{figure sim res A1}
\end{figure}
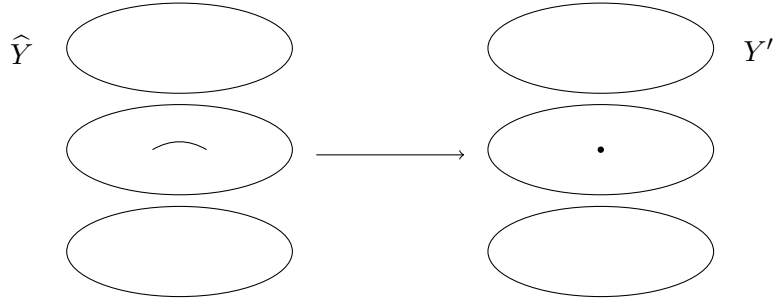

\subsection{Other cases}
In the $A_1$ case above, the resolution introduced a single complete curve over $0\in U$, and was an isomorphism over $U \backslash \{0\}$. In general, the locus in $U$ parametrising those fibres in which complete curves are introduced is more complicated. It is known as the \emph{discriminant locus}.

\begin{example}\label{example A2 versal} The discriminant locus for the $A_2$ case is illustrated in Figure~\ref{figure versal A2}, where $U = \{t \in \mathbb{C}^3 : t_1 + t_2 + t_3 = 0\}, $ with the natural action of $W = S_3$. A~selection of fibres of the simultaneous resolution $\simRes$ over $U$ are shown. The resolution $\surfaceRes$ is the fibre over the origin in $U$. The two intersecting projective lines in the exceptional locus of $\surfaceRes$ deform over two components of the discriminant locus~$\discriminantLocus$, and the two lines deform to a single smooth curve over a third component.
\end{example}

\newcommand{\defDeformationPos}{4}
\newcommand{\defParamSpaceSize}{4}
\newcommand{\defParamSpaceStretch}{1.5}
\newcommand{\defLabelSep}{-0.1}
\begin{figure}[hb]
\begin{center}
\begin{tikzpicture}[scale=0.7]
\node (A) at (0,0,0) {\begin{tikzpicture}[scale=\defSurfaceScale]\defTypeASurface{2}{}\end{tikzpicture}};
\node (Aleft) at (0,\defDeformationPos,0) {\begin{tikzpicture}[scale=\defSurfaceScale]\defTypeASurfaceShifted{1}{-1}\end{tikzpicture}};
\node (Aright) at (0,0,-\defDeformationPos*\defParamSpaceStretch) {\begin{tikzpicture}[scale=\defSurfaceScale]\defTypeASurfaceShifted{1}{1}\end{tikzpicture}};
\node (Amiddle) at (0,-\defDeformationPos,-\defDeformationPos*\defParamSpaceStretch) {\begin{tikzpicture}[scale=\defSurfaceScale]\defTypeASurface{1}{}\end{tikzpicture}};
\node (Agen) at (0,-\defDeformationPos/2,\defDeformationPos*\defParamSpaceStretch/2) {\begin{tikzpicture}[scale=\defSurfaceScale]\defSurface{3}{}\end{tikzpicture}};

\draw[gray](8,\defParamSpaceSize,\defParamSpaceSize*\defParamSpaceStretch) -- (8,-\defParamSpaceSize,\defParamSpaceSize*\defParamSpaceStretch) -- (8,-\defParamSpaceSize,-\defParamSpaceSize*\defParamSpaceStretch) -- (8,\defParamSpaceSize,-\defParamSpaceSize*\defParamSpaceStretch) -- (8,\defParamSpaceSize,\defParamSpaceSize*\defParamSpaceStretch); 
\draw(8,\defParamSpaceSize,0) -- (8,-\defParamSpaceSize,0); 
\draw(8,0,\defParamSpaceSize*\defParamSpaceStretch) -- (8,0,-\defParamSpaceSize*\defParamSpaceStretch); 
\draw(8,\defParamSpaceSize,\defParamSpaceSize*\defParamSpaceStretch) -- (8,-\defParamSpaceSize,-\defParamSpaceSize*\defParamSpaceStretch); 
	
\node (P) at (8,0,0) [Point] {};
\node [label={[shift={(0,-0.025)}]\scriptsize $t_1=t_2$}] (Pleft) at (8,\defDeformationPos,0) [Point] {};
\node [label={[shift={(0.6,-0.275)}]\scriptsize $t_2=t_3$}] (Pright) at (8,0,-\defDeformationPos*\defParamSpaceStretch) [Point] {};
\node [label={[shift={(0.6,-0.275)}]\scriptsize $t_1=t_3$}] (Pmiddle) at (8,-\defDeformationPos,-\defDeformationPos*\defParamSpaceStretch) [Point] {};
\node (Pgen) at (8,-\defDeformationPos/2,\defDeformationPos*\defParamSpaceStretch/2) [Point] {};

\node[right] (U) at (8,\defDeformationPos,-\defDeformationPos*\defParamSpaceStretch) {$U$};
\node[shift={(-0.2,0)},below right] (L) at (8,0,-\defDeformationPos*\defParamSpaceStretch/2) {$L$};

\draw[-,gray,dashed] (A) -- (P);
\draw[-,gray,dashed] (Aleft) -- (Pleft);
\draw[-,gray,dashed] (Aright) -- (Pright);
\draw[-,gray,dashed] (Amiddle) -- (Pmiddle);
\draw[-,gray,dashed] (Agen) -- (Pgen);

\node (labelA) at (A.north west) {$\surfaceRes$};
\node[shift={(-1.6,-0.3)},above of=Aleft] (labelZ) at (Aleft.north) {$\simRes$};

\end{tikzpicture}
\end{center}
\caption{Simultaneous resolution $\simRes$ for type $A_2$}
\label{figure versal A2}
\end{figure}
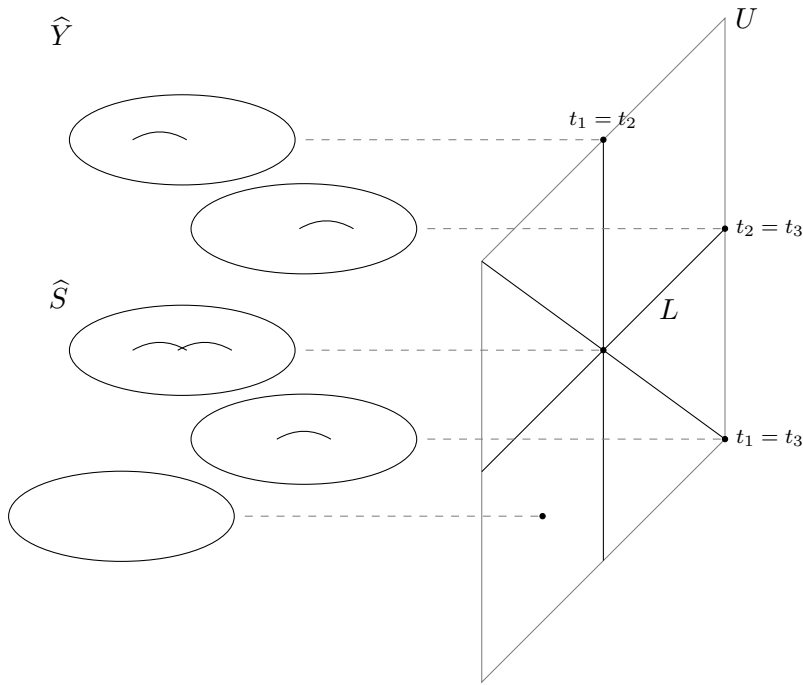

\subsection{Partial resolutions}\label{section partial res}

The resolution $\surfaceRes$ of Section~\ref{du val} is minimal amongst the resolutions of the singular surface $S$. There also exist partial resolutions $\surface'$ of~$\surface$, which may be constructed by contracting choices of curves in $\surfaceRes$ corresponding to vertices~$\chosenVertices \subseteq \DynkinDiag$ in the Dynkin diagram. The deformations of these partial resolutions have a description as follows \cite{KM}. Take $W$ the Weyl group acting on~$U$, and now let $\WeylSubgroup$ be the subgroup of $W$ associated to the vertices~$\chosenVertices$. Taking~$Y$ once again to be the versal deformation of the singularity $S$ over the base $U/W$ and now base changing to $U/W'$, there exists a partial resolution of the singularities of~$Y$, yielding a versal deformation of $S'$. Denoting this by $Y'$, there is a commutative diagram of Cartesian squares as follows, and the partial resolution $S' \to S$ is recovered by taking the fibre of the top line over~$0 \in U/W$.
\[
\begin{tikzcd}
\simRes \ar{r} \ar{d} & Y' \ar{r} \ar{d} & Y \ar{d} \\
U \ar{r} & U/\WeylSubgroup \ar{r} & U/W
\end{tikzcd}
\]

\begin{figure}[h]
\begin{center}
\begin{tikzpicture}[scale=0.7]

\node (D) at (0,0) {\begin{tikzpicture}[scale=\defSurfaceScale]\defTypeDSurface{3}\end{tikzpicture}};
\node (Dpar) at (7,0) {\begin{tikzpicture}[scale=\defSurfaceScale]\defTypeASurface{1}{}\end{tikzpicture}};
\node (Dcon) at (14,0) {\begin{tikzpicture}[scale=\defSurfaceScale]\defSingSurface{3}{}\end{tikzpicture}};
\node (DDynkin) at (0,0+\defDynkinSep) {\begin{tikzpicture}[scale=\defDynkinScale]
\defTypeDDynkin{3}\end{tikzpicture}};
\node (DDynkinpar) at (7,0+\defDynkinSep) {\begin{tikzpicture}[scale=\defDynkinScale]
\defTypeDDynkinSingle{3}
\end{tikzpicture}};

\draw[->] (D.east) -- (Dpar);
\draw[->] (Dpar) -- (Dcon);

\node (labelD) at (D.north west) {$\surfaceRes$};
\node (labelDpar) at (Dpar.north west) {$\surface'$};
\node (labelDcon) at (Dcon.north west) {$\surface$};

\end{tikzpicture}
\end{center}
\caption{Partial resolution for type $D_4$}
\label{figure partial resolution D4}
\end{figure}
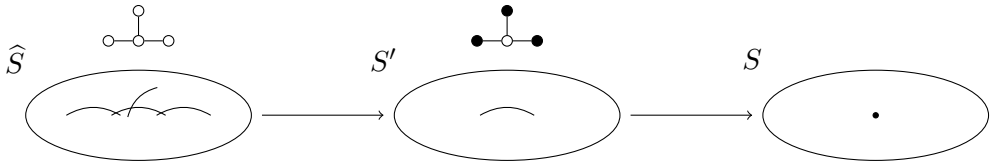

\begin{example}In the $D_4$ case illustrated in Figure~\ref{figure partial resolution D4}, we blow down the black curves $\chosenVertices$ in the Dynkin diagram $\DynkinDiag$ to leave a single remaining curve. This blowdown then factors the resolution $\surfaceRes \to \surface$, as shown.
\end{example}

\begin{example}\label{example A2 partial res} For the $A_2$ case in Example~\ref{example A2 versal}, choosing $\DynkinDiag'$ to be one of the two curves, $Y'$ may be constructed by contracting the codimension~$2$ locus in $\simRes$ consisting of a deformation of that curve lying over the corresponding component of~$L$, and then quotienting by $W' = S_2$.\end{example}

As before, a generic fibre of $Y'$ contains no complete curves, but $U/W'$ has a locus parametrising those deformations which do contain such curves, or singularities. This locus is given by the quotient of the discriminant locus $\discriminantLocus$ in $U$, and denoted also by $\discriminantLocus$.

\section{Contractions of $3$-folds}
\label{section threefolds}

We now return to considering a $3$-fold contraction of a curve $C \cong \mathbb{P}^1$, as in the introduction. I review the description of the geometry of such curves, and deformations of their contractions, which may be used to replace a non-generic curve $C$ with a collection of generic curves, following \cite{Pagoda} and \cite{Bryan-Katz-Leung}.

\subsection{Setup}
\label{section setup}
A \emph{contraction} is here taken to be a projective birational morphism $\contraction\colon X \to X_{\con}$ of $3$-folds, with the singularities of $X_{\con}$ being rational. A \emph{contraction of a curve} $C$ is a contraction which maps a curve $C$ to a point~$p$, and is an isomorphism from $X \backslash C$ onto $X_{\con} \!\backslash p$.
We say that $\contraction$ is a \emph{flopping~contraction} if the canonical bundle $K_X$ is trivial in a neighbourhood of~$C$.
In the following, $C \cong \mathbb{P}^1$.

For simplicity, I take $X$ to be smooth and $X_{\con} = \Spec R$, a~complete local ring, though the theorems presented later also work in a global setting, and under weaker assumptions on singularities, which I identify in Remarks~\ref{remark contraction algebra assumptions} and~\ref{remark contraction theorem assumptions}.

\begin{remark} In the situation above, the scheme-theoretic fibre $f^{-1}(p)$ may be a finite-order thickening of $C$. This is key to much of the geometric subtlety which follows.\end{remark}

\subsection{Slicing the singularity}
Given a flopping contraction $\contraction$ of a curve $C \cong \mathbb{P}^1$ to a point $p$, we take a generic hypersurface $\hypersurface$ through $p \in X_{\con}$.
A result of Reid \cite[Theorems~1.1,~1.14]{Pagoda} gives that $\hypersurface$ has a Kleinian singularity at $p$, so that we may apply the discussion of Section~\ref{section surfaces}.
Write $\hypersurfacePRes$ for the preimage of $H$ under $\contraction$, so that there is a map 
\begin{equation*}
\hypersurfaceResMap \colon \hypersurfacePRes \longrightarrow \hypersurface \ni p.
\end{equation*}
Now the reduced fibre of $\hypersurfaceResMap$ above $p$ is exactly $C \cong \mathbb{P}^1$, and so this map cannot be the resolution $\hypersurfaceRes \to \hypersurface$, unless the singularity of $\hypersurface$ is of the simplest Dynkin type~$A_1$. Indeed in general it is a partial resolution, factoring $\hypersurfaceRes \to \hypersurface$ as follows:
\begin{equation*}
\hypersurfaceRes  \overset{q}{\longrightarrow} \hypersurfacePRes \overset{\hypersurfaceResMap}{\longrightarrow} \hypersurface.
\end{equation*}
Here $q$ is a blowdown of certain projective lines in $\hypersurfaceRes$ as in Section~\ref{section partial res}, and identifies one of these lines with our curve $C$.

Recalling that the projective lines in $\hypersurfaceRes$ correspond to the vertices of a Dynkin diagram, we thereby associate to our contraction $\contraction$ the data of a Dynkin diagram $\DynkinDiag$ with a marked point $i$. This data is referred to as the \emph{Dynkin type} of the floppable curve~$C$.

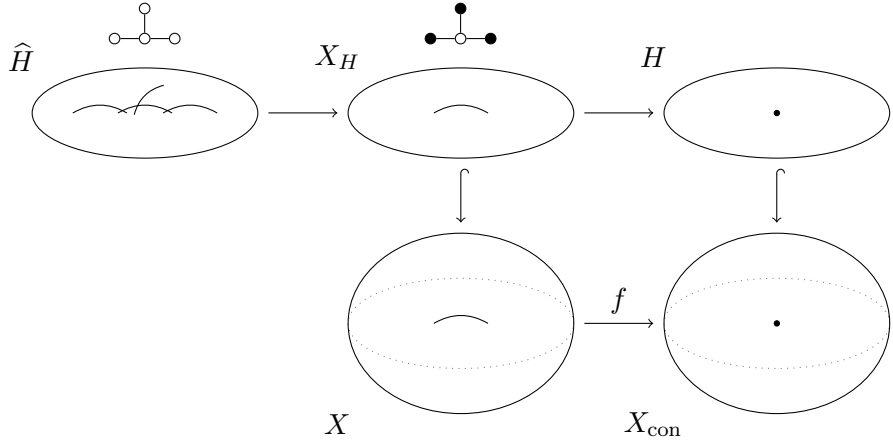
\begin{figure}[ht]
\begin{center}
\begin{tikzpicture}[scale=0.7]
\node (D) at (0,4) {\begin{tikzpicture}[scale=\defSurfaceScale]\defTypeDSurface{3}\end{tikzpicture}};
\node (Dpar) at (6,4) {\begin{tikzpicture}[scale=\defSurfaceScale]\defTypeASurface{1}{}\end{tikzpicture}};
\node (Dcon) at (12,4) {\begin{tikzpicture}[scale=\defSurfaceScale]\defSingSurface{3}{}\end{tikzpicture}};

\node (DDynkin) at (0,4+\defDynkinSep) {\begin{tikzpicture}[scale=\defDynkinScale]
\defTypeDDynkin{3}\end{tikzpicture}};
\node (DDynkinpar) at (6,4+\defDynkinSep) {\begin{tikzpicture}[scale=\defDynkinScale]
\defTypeDDynkinSingle{3}\end{tikzpicture}};

\node (ThreefoldDpar) at (6,0) {\begin{tikzpicture}[scale=\defSurfaceScale]\defTypeASurface{1}{gray, dotted}\defThreefold{3}\end{tikzpicture}};
\node (ThreefoldDcon) at (12,0) {\begin{tikzpicture}[scale=\defSurfaceScale]\defSingSurface{3}{gray, dotted}\defThreefold{3}\end{tikzpicture}};

\draw[->] (D) -- (Dpar);
\draw[->] (Dpar) -- (Dcon);
\draw[right hook->] (Dpar) -- (ThreefoldDpar);
\draw[right hook->] (Dcon) -- (ThreefoldDcon);
\draw[->] (ThreefoldDpar) -- node[above] {$\contraction$} (ThreefoldDcon);

\node (labelD) at (D.north west) {$\hypersurfaceRes$};
\node (labelDpar) at (Dpar.north west) {$\hypersurfacePRes$};
\node (labelDcon) at (Dcon.north west) {$\hypersurface$};
\node (labelThreefoldDpar) at (ThreefoldDpar.south west) {$\threefold$};
\node (labelThreefoldDcon) at (ThreefoldDcon.south west) {$\threefoldCon$};
\end{tikzpicture}
\end{center}
\caption{Slicing a singularity by a generic hyperplane, type $D_4$}
\label{figure generic hyperplane D4}
\end{figure}

\begin{example}\label{example type D}\cite[3.14]{DW1} I present a family of examples of type $D_4$ flopping curves, which we will return to throughout. For each $k \in \mathbb{N}$, put 
\[
X_{k,\con} = \big\{ u^2+v^2y=x(x^2+y^{2k+1}) \big\} \subset \mathbb{C}^4.
\]
This singularity can be blown up in an ideal to give a threefold $X_k$ with a single complete curve $C \cong \mathbb{P}^1$ \cite[Example 5.15]{Pagoda}. Figure~\ref{figure generic hyperplane D4} depicts the generic hypersurface section~$\hypersurface$ through the singularity, along with its partial resolution~$\hypersurfacePRes$, and full resolution~$\hypersurfaceRes$.
\end{example}

It turns out that only flopping curves of types $A_1$, $D_4$, $E_6$, $E_7$, and $E_8$ may arise. Furthermore, in the cases $D_4$, $E_6$, and $E_7$, the marked point must be the central vertex of the Dynkin diagram, whereas in the case $E_8$, the marked point may move one place along the long arm of the diagram. See Katz--Morrison \cite{KM} or Kawamata \cite{Kawa} for an account of the possibilities.

\subsection{Characterising the $3$-fold}
We now obtain the desired description of the contraction $\contraction$ in terms of deformations of the surfaces in Section~\ref{section surfaces}, as follows.

Writing $\hypersurfaceEq=0$ for the equation of our hypersurface $\hypersurface \subset X_{\con}$, the $3$-fold $\threefold$ may be viewed as a family over a formal neighbourhood $\neighbourhood \subset \mathbb{A}^1$ via $\hypersurfaceEq \circ \contraction$. Now indeed this is a deformation of a partial resolution $S'$ associated to the marked Dynkin diagram $(\DynkinDiag, i)$, and so by the description of the versal deformation $Y'$ in Section~\ref{section partial res} we obtain a classifying map $g \colon \neighbourhood \to U/\WeylSubgroup$ for an appropriate subgroup $W'$ of the Weyl~group. This may be viewed as a germ of a $\mathbb{C}$-parametrised path in the orbit space~$U/W'$. 

To summarise, we have associated to a flopping contraction of a curve $C \cong \mathbb{P}^1$, the data $(\DynkinDiag, i, g)$ of a Dynkin diagram with a marked vertex, and a path of $W'$-orbits in the associated $U$, as follows. \begin{equation*}(X,f) \quad \longleftrightarrow \quad (\DynkinDiag, i, g)\end{equation*} In the other direction, the data~$(\DynkinDiag, i, g)$ may be used to construct a contraction, by pulling back the versal deformation $Y'$ via $g$.

\begin{remark}For the $D_4$ case of Example~\ref{example type D}, the classifying map to the versal deformation base is given explicitly by Aspinwall and Morrison~\cite[Section~4.2]{AM}.
\end{remark}

\subsection{Case $A_1$}
\label{section type A}

I spell this case out as it gives the easiest class of examples. Here $\DynkinDiag$ consists of a single vertex, and we must therefore take this for the marked vertex~$i$. Then $W' = \{ \Id \}$ with $U = \mathbb{C}$. Our path is therefore given by a function $g\colon \mathbb{C} \to U = \mathbb{C}$. Taking $g(\pathParameter)=\pathParameter^d$, base change of \eqref{equation conifold} in our previous type $A_1$ example gives the equation of a singularity \begin{equation}\label{equation width d singularity}\big\{x^2 + y^2 = z^2 - \pathParameter^{2d}\big\} \subset \mathbb{C}^4.\end{equation}

For $d=1$, the $3$-fold $X$ is just the simultaneous resolution of \eqref{equation conifold} discussed in Section~\ref{section A1 case}, containing a single $(-1,-1)$-curve. For $d\geq 2$, base change gives a simultaneous resolution of the family of singularities \eqref{equation width d singularity} with a single complete curve $C$ in the fibre $w=0$. This curve has normal bundle $\cO_C(-2) \oplus \cO_C$, and is referred to as a $(-2,0)$-curve. The number $d$ is called the \emph{width} by Reid  \cite[Definition~5.3]{Pagoda}.

\subsection{Deforming the $3$-fold}\label{section deforming 3-fold} We may now take $\tilde{g}$, a new path in $U/W'$ transverse to the discriminant locus $L$, whose germ is a small deformation of $g$, as in Figure~\ref{figure general classifying map}. Base changing the versal deformation $Y'$  by $\tilde{g}$, a new deformed $3$-fold $\threefoldDef$ is obtained. Similarly, a singular base $\threefoldDefCon$ is obtained by base changing the versal deformation~$Y$.

\newcommand{\DefAngle}{70}
\newcommand{\DefDeformationAngle}{40}
\newcommand{\DefDeformationShift}{0.5}
\newcommand{\DefSurfaceShrink}{0.8}
\begin{figure}[h]
\begin{center}
\begin{tikzpicture}[scale=1.5]

\node (ClassifyingMap) at (0,0) {
\begin{tikzpicture}[scale=1.5,rotate=-60]
\draw[thick] (0,0) -- (1.2,0);
\node[left] at (1.2,0) {$L$};

\draw [thick] (-1,-1) to[out=\DefAngle,in=180] (0,0) to[out=0,in=180+\DefAngle] (1,1);
\node[right] at (1,1) {$g$};

\draw [thick] (-1,-1+\DefDeformationShift) to[out=\DefAngle,in=180-\DefDeformationAngle/2] (0,0) to[out=-\DefDeformationAngle/2,in=180+\DefAngle] (1,1-\DefDeformationShift);
\node[right] at (1,1-\DefDeformationShift) {$\tilde{g}$};
\end{tikzpicture}
};	

\end{tikzpicture}
\end{center}
\caption{Classifying map $g$ and deformation $\tilde{g}$}
\label{figure general classifying map}
\end{figure}
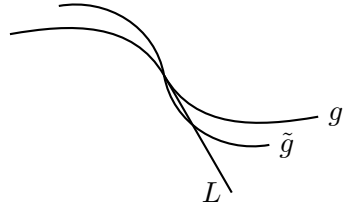

\begin{figure}[h]
\begin{center}
\begin{tikzpicture}[scale=0.7]

\node (ClassifyingMap) at (0,0) {
\begin{tikzpicture}[scale=1.4]
\draw[->] (-1,0) -- (1,0);
\draw[->] (0,-1) -- (0,1);
\draw [thick] (-1,-1) to[out=\DefAngle,in=180] (0,0) to[out=0,in=180+\DefAngle] (1,1);
\node[right] at (1,0) {\scriptsize$\pathParameter$};
\node[above] at (0,1) {\scriptsize$g(\pathParameter)$};
\end{tikzpicture}
};	

\node (OrigThreefold) at (4.5,0) {$\begin{array}{c}
\begin{tikzpicture}[scale=\defSurfaceScale*\DefSurfaceShrink]\defSurface{3}{}\end{tikzpicture} \\
\begin{tikzpicture}[scale=\defSurfaceScale*\DefSurfaceShrink]\defTypeASurface{1}{}\end{tikzpicture} \\
\begin{tikzpicture}[scale=\defSurfaceScale*\DefSurfaceShrink]\defSurface{3}{}\end{tikzpicture}
\end{array}$};

\node (DefClassifyingMap) at (10,0) {
\begin{tikzpicture}[scale=1.4]
\draw[->] (-1,0) -- (1,0);
\draw[->] (0,-1) -- (0,1);
\draw [thick] (-1,-1+\DefDeformationShift) to[out=\DefAngle,in=180-\DefDeformationAngle/2] (0,0) to[out=-\DefDeformationAngle/2,in=180+\DefAngle] (1,1-\DefDeformationShift);
\node[right] at (1,0) {\scriptsize$\pathParameter$};
\node[above] at (0,1) {\scriptsize$\tilde{g}(\pathParameter)$};
\end{tikzpicture}
};	

\node (DefThreefold) at (14.5,0) {$\begin{array}{c}
\begin{tikzpicture}[scale=\defSurfaceScale*\DefSurfaceShrink]\defTypeASurface{1}{}\end{tikzpicture} \\
\begin{tikzpicture}[scale=\defSurfaceScale*\DefSurfaceShrink]\defTypeASurface{1}{}\end{tikzpicture} \\
\begin{tikzpicture}[scale=\defSurfaceScale*\DefSurfaceShrink]\defTypeASurface{1}{}\end{tikzpicture}
\end{array}$};

\node (labelOrigThreefold) at (2,-3) {\small $(-2,0)$-curve, width $3$};
\node (labelDefThreefold) at (12,-3) {\small $3$ disjoint $(-1,-1)$-curves};
\node (labelOrig) at (4.5,3) {\small $X$};
\node (labelDef) at (14.5,3) {\small $\tilde{X}$};

\end{tikzpicture}
\end{center}
\caption{Deformation of classifying map $g$ for $(-2,0)$-curve}
\label{figure type A classifying map}
\end{figure}
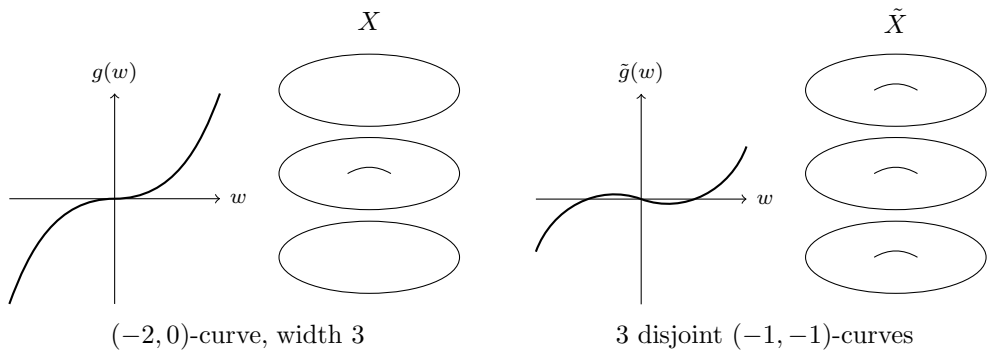

In the $A_1$ case of the previous Section~\ref{section type A} the discriminant locus $L=\{0\} \subset \mathbb{C}$. A non-transverse classifying map $g$ with $d>1$ may be deformed to a transverse map~$\tilde{g}$. The roots of $\tilde{g}$ then correspond to $d$ disjoint $(-1,-1)$-curves in the deformed $3$-fold~$\threefoldDef$. This is sketched in Figure~\ref{figure type A classifying map}.

Generally, non-transverse intersections with the discriminant locus in $U/\WeylSubgroup$ correspond to non-generic curves, whereas transverse intersections correspond to $(-1,-1)$-curves. I now indicate how such intersections may be counted.

\subsection{Counting curves in the deformation}\label{section curves deformation}
The curves in the deformed $3$-fold~$\tilde{X}$ may be characterised in terms of certain enumerative invariants on $X$. As preparation, we produce a sequence of infinitesimal neighbourhoods of $C \subset \threefold$ as follows. Taking $\mathcal{I}$ to be the ideal of definition of $C \subset \hypersurfacePRes$, we define a thickened curve~$C^{(j)}$ for $j \in \mathbb{N}$ by the $j^\text{th}$ symbolic power $\mathcal{I}^{(j)}$ of $\mathcal{I}$, viewed as a subscheme of $\threefold$.

Each $C^{(j)}$ corresponds to an isolated point in $\operatorname{Hilb}(\threefold)$, the Hilbert scheme of closed subschemes of $\threefold$. Write $n_j$ for the multiplicity of this point. These~$n_j$ coincide with the genus zero \GV invariants associated to~$C$, as defined by Katz~\cite{Katz}.

Using the multiplicities $n_j$, Bryan--Katz--Leung \cite{Bryan-Katz-Leung} describe intersections of a transverse deformation $\tilde{g}$ with the components of the discriminant locus~$\discriminantLocus$ in~$U/W'$, and also show that the curves in the associated $3$-fold $\threefoldDef$ are all $(-1,-1)$-curves, with different curve classes $j[C]$ corresponding to intersections of $\tilde{g}$ with different components of~$L$. This is achieved by a careful analysis of the geometry of the discriminant locus in different Dynkin types, using explicit presentations of the simultaneous resolutions due to Katz--Morrison \cite{KM}. Their result is the following.

\begin{proposition}\cite[\S 2.1]{Bryan-Katz-Leung}\label{proposition BKL intersection} Writing $n_j$ for the multiplicities associated to the curve $C \subset X$ as above, the exceptional locus of a generic deformed contraction $\threefoldDef \to \threefoldDefCon$ consists of $(-1,-1)$-curves, with $n_1$ curves of class~$[C]$, $n_2$ of class~$2[C]$, and so on, up to~$l[C]$.
\end{proposition}

Here the \emph{length} $l$ of $C$ is the length of the scheme-theoretic fibre $f^{-1}(p)$ at the generic point of $C$. In type $A$, the length $l=1$, and the multiplicity $n_1$ coincides with the width $d$ of $C$. In types $D_4$, $E_6$, and $E_7$, the length $l=2$, $3$, and $4$, respectively. In type $E_8$ the length takes values~$5$ or $6,$ if the marked point of the Dynkin diagram is adjacent to the central vertex (on the long arm), or the central vertex itself, respectively.

\section{The contraction algebra $\CA$}
\label{section contraction algebra}

Having explored how isolated contractible curves arise, we introduce a new invariant associated to them. As an example application we will calculate the \GV invariants $n_j$ for our $D_4$ example in the following Section~\ref{section numerical invariants}, using a result of Toda \cite{TodaGV}.

As in the previous Section~\ref{section threefolds}, we consider a $3$-fold contraction~$\contraction\colon X \to X_{\con}=\operatorname{Spec} R$, a complete local ring, contracting a curve~$C \cong \mathbb{P}^1$. Unlike the previous section, the following definitions and theorems do not require the contraction~$\contraction$ to be flopping: indeed they extend to the case of singular flops and flips (see Remarks~\ref{remark contraction algebra assumptions} and~\ref{remark contraction theorem assumptions}).
\subsection{Via a tilting algebra}

In this setting there exists a tilting bundle~$\mathcal{T}$ on $\threefold$ by a construction of Van~den~Bergh \cite[Theorem~A]{VdB1d} given explicitly in \cite[\S 2.3]{DW1}. The usefulness of this bundle stems from the fact that the algebra $\AB := \End_X(\mathcal{T})$ is derived equivalent to $\threefold$. The bundle $\mathcal{T}$ decomposes as $\mathcal{T} = \cO_X \oplus \mathcal{N}$. Writing $N$ for the $R$-module of global sections of~$\mathcal{N}$, we find that $\AB\cong\End_R(R \oplus N)$ using~\cite[4.2.1]{VdB1d}. The \emph{contraction~algebra}~$\CA$ may then be defined as follows.

\begin{definition}\cite[\S 2.3--\S 2.4]{DW1}\label{Definition contraction algebra} For a contraction of a curve $C \cong \mathbb{P}^1$ as in Section~\ref{section setup}, put $\CA = \AB / I_\con$, where $I_\con$ is the two-sided ideal of~$\AB$ generated by endomorphisms $\upphi \in \AB\cong\End_R(R \oplus N)$ which factor through a direct summand of a free $R$-module. \end{definition}

\begin{remark}\label{remark contraction algebra assumptions}We do not need smoothness of $X$ to make this definition, or for the Theorem~\ref{theorem representability} which follows. These are treated in \cite{DW1} under the weaker assumption that $X$ is normal, with only Cohen--Macaulay canonical singularities.
\end{remark} 

Examples of contraction algebras $\CA$ are given in Section~\ref{section Acon examples}.

\subsection{Via noncommutative deformations}

The contraction algebra is also characterised, up to isomorphism, by the noncommutative deformations of the curve~$C$. Roughly speaking, $\CA$ plays the role of an algebra of functions on a noncommutative space parametrising the deformations of the sheaf $\cO_C$ in coherent sheaves on~$\threefold$. To be precise, we find that $\CA$ represents a certain deformation functor. I give a partial statement of the required formalism here, referring the reader to \cite[\S 2]{DW1} for full details.

First we take a category $\art_1$ of finite-dimensional $\mathbb{C}$-algebras $\TestAlg$, not necessarily commutative. These are endowed with extra structure, in particular an ideal~$\mathfrak{n}$ which should be considered as the unique closed point of the noncommutative space corresponding to $\TestAlg$. For a given coherent sheaf $E$ on $\threefold$ we then construct a \emph{noncommutative deformation functor} 
\[
\Def_E\colon \art_1 \to \Sets,
\]
taking $(\TestAlg,\mathfrak{n})$ to the set of $\TestAlg$-families $\cF$ which recover $E$ after reduction modulo the ideal $\mathfrak{n}$. As part of their definition, such $\cF$ are objects of $\coh(\cO_X \otimes_{\mathbb{C}} \TestAlg)$.

Note that restricting $\Def_E$ to a category of commutative $\mathbb{C}$-algebras, the usual deformation functor for the sheaf $E$ is obtained. This is referred to here as the \emph{classical deformation functor}, and denoted by $\cDef_E$. We then have the following.

\begin{theorem}\cite[1.1, 2.13(1), \S 3.1]{DW1}\label{theorem representability} Given a contraction of a curve $C \cong \mathbb{P}^1$ as in Section~\ref{section threefolds}, the noncommutative deformation functor $\Def_{\cO_C(-1)}$ is represented by the contraction algebra $\CA$ of Definition~\ref{Definition contraction algebra}, and this algebra is finite-dimensional. Furthermore, the classical deformation functor $\cDef_{\cO_C(-1)}$ is represented by the abelianisation $\CAab$, given by quotienting $\CA$ by its two-sided ideal of commutators.
\end{theorem}

\begin{remark}\label{remark contraction theorem assumptions}The statement above is a special case of the theorem in \cite{DW1}. There we only ask that $X$ is normal and Cohen--Macaulay canonical, and we allow $X_{\con}$ to be quasi-projective, rather than complete local. In particular, the theorem covers both flops and flips in singular settings.
\end{remark}

\subsection{Results and examples}\label{section Acon examples}

The following theorem may be obtained using knowledge of the number of generators and relations of the algebra $\CA$, and applying a powerful result of Smoktunowicz \cite{Agata}, or alternatively by a careful analysis of hyperplane sections. Both approaches are discussed in \cite[Section~3.4]{DW1}. 

\begin{theorem}\cite[3.15]{DW1}\label{Acon commutative} The contraction algebra $\CA$ is commutative if and only if the contracted curve $C$ is of Dynkin type $A$.\end{theorem}

\begin{example}\cite[3.10]{DW1} If a flopping curve $C$ is of type $A$ (as in Section~\ref{section type A} for example) then $\CA\cong \mathbb{C}[\upepsilon]/\upepsilon^d$, where $d$ is the width of $C$.
\end{example}

\begin{example}\cite[3.14]{DW1}\label{example contraction algebra D4} For the $D_4$ case of Example~\ref{example type D}, 
\[
\CA \cong \frac{\mathbb{C}\langle x,y\rangle}{xy=-yx\mbox{, }x^2=y^{2k+1}}.
\] 
A vector space basis for $\CA$ is given by the monomials $x^a y^b$ with $0 \leq a \leq 2$ and $0 \leq b \leq 2k$. We have
\[
\CAab \cong \frac{\mathbb{C}[x,y]}{xy=0\mbox{, }x^2=y^{2k+1}}.
\] 
Both these algebras, and the associated deformations, are illustrated in Figure~\ref{figure noncomm deformations}.
\end{example}

\def\noncommShift{\size*6}
\def\fibreShift{-1.5}
\begin{figure}[h]
\centering
\begin{tikzpicture}[scale=2,node distance=6]
\node at (0,-0.6) {$
  \begin{tikzpicture}[xscale=0.3,yscale=0.6]
    \defThreeFold{\size*3}{0}{black};
    \defCurve{\size*3}{0}{black};
    \node at (\size*2,\size/4) {$\small X$};
    \defGridLines{\fibreShift}{0}{5}{3};
    \defUniversalComm{\fibreShift}{0}{2}{4};
    \draw[->] (0.75*\size,0) -- (\size*1.75,0);

    \draw[->] (\fibreShift+0.25,\gridShift+0.25,0) -- node[above,pos=0.5] {$\scriptstyle \varepsilon$} (\fibreShift+1,\gridShift+0.25,0);
\node at (\fibreShift-2.5,-4.5) {$\small \CAab$};
\end{tikzpicture}$};

\node at (\size*0.68,-0.6) {$
  \begin{tikzpicture}[xscale=0.3,yscale=0.6]
    \defThreeFold{\noncommShift}{0}{white};
    \defGridLines{\noncommShift}{0}{5}{3};
    \defUniversal{\noncommShift}{0}{2}{4};
    \draw[->] (\noncommShift-0.5*\size,0) -- (\noncommShift-\size*1.5,0);
    \draw[->] (\noncommShift-0.6,\gridShift,0) -- node[left,pos=0.1] {$\scriptstyle x$} (\noncommShift-0.6,\gridShift,1);
    \draw[->] (\noncommShift+0.25,\gridShift+0.25,0) -- node[above,pos=0.5] {$\scriptstyle y$} (\noncommShift+1,\gridShift+0.25,0);
\node at (\noncommShift+6.5,-4.5) {$\small \CA$};
\end{tikzpicture}$};

\end{tikzpicture} 
\caption{Deformations of type $D_4$ curve for $k=2$}
\label{figure noncomm deformations}
\end{figure}
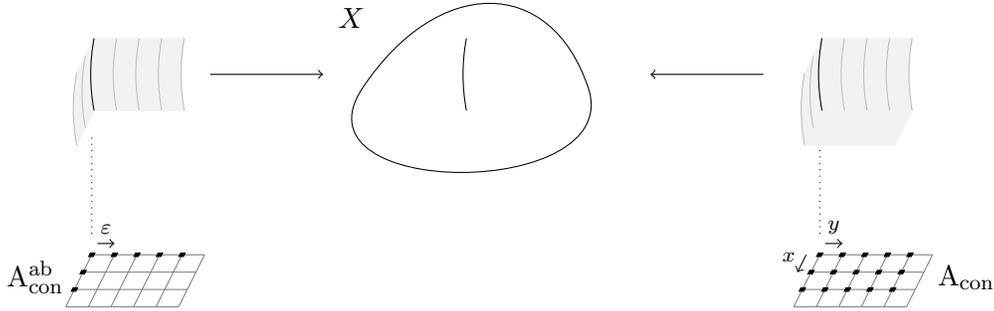

\section{Numerical invariants}
\label{section numerical invariants}

\subsection{Widths} We refer to the dimensions of the contraction algebra $\CA$ and its abelianisation $\CAab$ as the \emph{noncommutative width} and \emph{commutative width} respectively, and denote these numbers by $\operatorname{wid}(C)$ and $\operatorname{cwid}(C).$  For the type~$A$ flopping curves of Section~\ref{section type A}, that is the case in which the width of Reid is defined, we find that all three notions of width coincide. More generally, the noncommutative and commutative widths may differ, as shown by the following example.

\begin{example}In our type $D_4$ example, $\operatorname{wid}(C) = 3(2k+1)$, and $\operatorname{cwid}(C) = 2k+3,$ from the results in Example~\ref{example contraction algebra D4}.
\end{example}

An analysis of hyperplane sections gives lower bounds for $\operatorname{wid}(C)$ according to Dynkin type, as explained in \cite[Remark~3.17]{DW1}. For type~$D_4$ the bound is~$4$, and for types~$E_6$, $E_7$, and $E_8$, the bounds are $12$, $24$, and $40$, respectively.

\begin{remark}The noncommutative width may also be obtained by a calculation of $\Ext$ groups on $X_{\con}$. See \cite[Remark~5.2]{DW1}.	
\end{remark}

\subsection{\GV invariants}
I now state the result of Toda from the introduction and explain how it sheds light on the width invariants defined above. As before, let $l$ be the length of the scheme-theoretic fibre $f^{-1}(p)$ at the generic point of $C$.
\begin{theorem}\cite{TodaGV}\label{theorem.toda} For a flopping contraction of a curve $C \cong \mathbb{P}^1$ as in Section~\ref{section threefolds}, on a $3$-fold $X$ with at worst Gorenstein terminal singularities, we have that
\[
\operatorname{wid}(C) = \sum_{j=1}^l j^2 n_j,
\]
and $\operatorname{cwid}(C) = n_1,$ where the $n_j$ are the multiplicities defined in Section~\ref{section curves deformation}, coinciding with genus zero \GV invariants.
\end{theorem}

Recalling from Proposition~\ref{proposition BKL intersection} that the $n_j$ also coincide with the numbers of $(-1,-1)$-curves of class~$j[C]$ in a generic deformation, we see that the noncommutative width~$\operatorname{wid}(C)$ counts the curves in such a deformation, with weights given by the squares of their curve classes.

\begin{example}For instance, in our running type $D_4$ case from Example~\ref{example contraction algebra D4}, we deduce that $n_1=2k+3$, and $n_2=k$. In particular, the multiplicity of the thickened curve $C^{(2)}$ in its Hilbert scheme is $k$. This is illustrated in Figure~\ref{fig.deformations2}. We also deduce that under a generic deformation, the flopping curve $C$ deforms to $2k+3$ curves of class $[C]$, and $k$ curves of class $2[C]$.
\end{example}

\def\noncommShift{\size*6}
\def\fibreShift{-1.5}
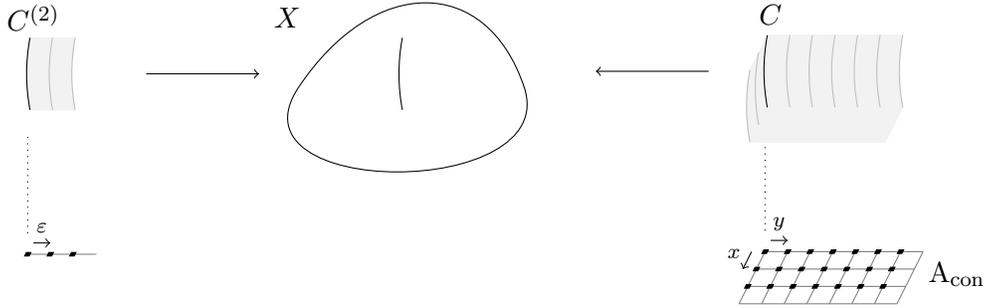
\begin{figure}[ht]
\centering
\begin{tikzpicture}[scale=2,node distance=6]
\node at (0,-0.45) {$
  \begin{tikzpicture}[xscale=0.3,yscale=0.6]
    \defThreeFold{\size*3}{0}{black};
    \defCurve{\size*3}{0}{black};
    \node at (\size*2,\size/4) {$\small X$};
    \defGridLines{\fibreShift}{0}{3}{0};
    \defUniversalComm{\fibreShift}{0}{0}{2};
    \draw[->] (0.75*\size,0) -- (\size*1.75,0);
    \draw[->] (\fibreShift+0.25,\gridShift+0.25,0) -- node[above,pos=0.5] {$\scriptstyle \varepsilon$} (\fibreShift+1,\gridShift+0.25,0);
\node at (\fibreShift+0.25,1.25) {$C^{(2)}$};
\end{tikzpicture}$};

\node at (\size*0.68,-0.6) {$
  \begin{tikzpicture}[xscale=0.3,yscale=0.6]
    \defThreeFold{\noncommShift}{0}{white};
    \defGridLines{\noncommShift}{0}{7}{3};
    \defUniversal{\noncommShift}{0}{2}{6};
    \draw[->] (\noncommShift-0.5*\size,0) -- (\noncommShift-\size*1.5,0);
    \draw[->] (\noncommShift-0.6,\gridShift,0) -- node[left,pos=0.1] {$\scriptstyle x$} (\noncommShift-0.6,\gridShift,1);
    \draw[->] (\noncommShift+0.25,\gridShift+0.25,0) -- node[above,pos=0.5] {$\scriptstyle y$} (\noncommShift+1,\gridShift+0.25,0);
\node at (\noncommShift+8.5,-4.5) {$\small \CA$};
\node at (\noncommShift+0.25,1.25) {$C$};
\end{tikzpicture}$};

\end{tikzpicture} 
\caption{Deformations of (thickened) type $D_4$ curve for $k=3$}
\label{fig.deformations2}
\end{figure}

\section{The derived category}
\label{section derived category}
In this section, I briefly indicate the application of $\CA$ to describing symmetries of the bounded derived category of coherent sheaves~$\D(X)$, and how this is used in Toda's proof of Theorem~\ref{theorem.toda}.

\subsection{Autoequivalences} Unpacking the statement of Theorem~\ref{theorem representability}, we obtain a \emph{universal family} $\cE \in \coh(\cO_X \otimes_{\mathbb{C}} \CA)$ for the noncommutative deformations of~$\cO_C(-1).$ Then we have the following theorem.

\begin{theorem}\label{theorem.twist}\cite[1.5, 7.14]{DW1} Given a flopping contraction of $C \cong \mathbb{P}^1$ for a projective $3$-fold~$X$ with at worst Gorenstein terminal singularities, there exists an autoequivalence $T_C$ of \,$\D(X)$ which acts on objects as follows.
 \[T_C(E) \cong \operatorname{Cone}\left( \mathbf{R}\!\Hom_X(\cE,E)\otimes_{\CA}^{\bf L} \cE \overset{\operatorname{ev}}{\longrightarrow} E \right)\]
\end{theorem}

\begin{remark}Replacing $\CA$ with the commutative algebra $\CAab$ in the above theorem does not give an autoequivalence of $\D(X)$, except in Type $A$ cases, where $\CA$ and $\CAab$ coincide. This is demonstrated in \cite[Theorem~3.15]{DW1}.\end{remark}

For a flop $X \birational X'$ there are equivalences of categories \[F\colon \D(X) \longleftrightarrow \D(X')\reversecolon F'\] given by the Fourier--Mukai transforms naturally associated to the graph of the birational map. This is a result of Bridgeland \cite{Bridgeland}, generalised by Chen to the Gorenstein terminal setting \cite{Chen}. We prove the following.

\begin{theorem}\label{theorem.twist-flop-flop}\cite[7.18]{DW1}
The autoequivalence $T_C$ of Theorem~\ref{theorem.twist} is an inverse to the flop--flop functor $F' \circ F$.\end{theorem}

\subsection{Application}

I finish with a brief description of Toda's proof \cite{TodaGV} of Theorem~\ref{theorem.toda}, relating the noncommutative width $\operatorname{wid}(C)$ and the invariants $n_j$. Interpolating the classifying maps $g\text{, }\tilde{g}\colon \mathbb{C} \to U/W'$ by a continuous family~$g^t$, we may obtain a $1$-parameter family of $3$-fold contractions $X^t \to X_{\con}^t$. (For this we should replace the versal deformations with suitable localisations, so as to recover the complete local fibre of the contraction that we started with, at $t=0$.) In this family $X^t$ there is a family flop of the exceptional locus, and thence a family of flop--flop functors. On~the central fibre, the flop--flop is inverse to $T_C$ by Theorem~\ref{theorem.twist-flop-flop}. On the generic fibre, the flop--flop is inverse to a composition of spherical twists \cite{ST}, one for each curve. Toda obtains the identities of Theorem~\ref{theorem.toda} by analysing the Fourier--Mukai kernels of these autoequivalences, and equating Hilbert polynomials.


\end{document}